\documentclass[12pt,final]{amsart}

\usepackage{graphicx} 
\usepackage[margin=1in]{geometry}
\usepackage{amsmath}

\usepackage{verbatim}
\usepackage{amssymb}
\usepackage{amscd}

\theoremstyle{definition}

\theoremstyle{remark}

\author{Tatyana Barron}
\address{Department of Mathematics, University of Western Ontario, London Ontario N6A 5B7, Canada}
\email{tatyana.barron@uwo.ca}

\author{Sarah Lanthier}
\address{Department of Biology, Department of Mathematics, University of Western Ontario, London Ontario N6A 5B7, Canada}
\email{slanthi2@uwo.ca}

\date{\today}

\title{Geometry and biological processes}
\begin{document}
\sloppy

\maketitle

\noindent {\bf Abstract.} We suggest a geometric approach to modeling biochemical processes, aiming at those processes that occur in humans with food sensitivities or chemical sensitivities.
  
  \
  

\

  
\

\section{Introduction} 
\begin{figure}[htb]
\begin{minipage}{0.32\textwidth}
\centering
\includegraphics[width=2.2in]{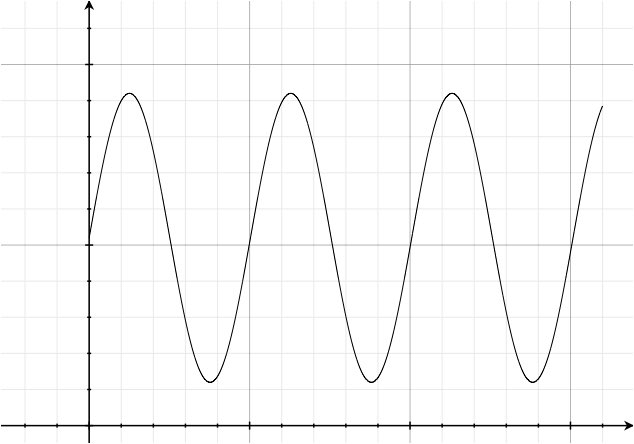}
\caption{} \label{figproc1}
\end{minipage}
\begin{minipage}{0.32\textwidth}
\centering
\includegraphics[width=2.2in]{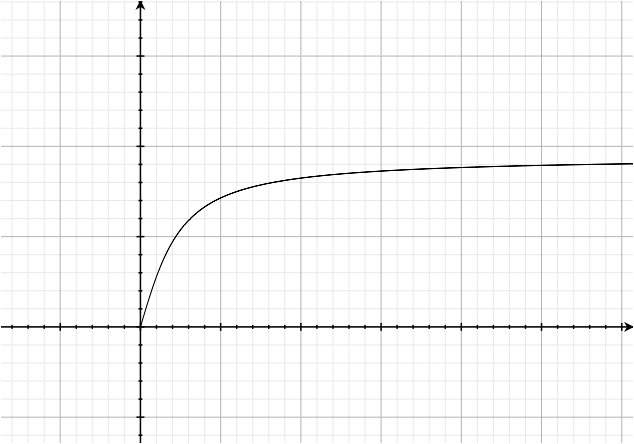}
\caption{} \label{figproc2}
\end{minipage}
\begin{minipage}{0.32\textwidth}
\centering
\includegraphics[width=2.2in]{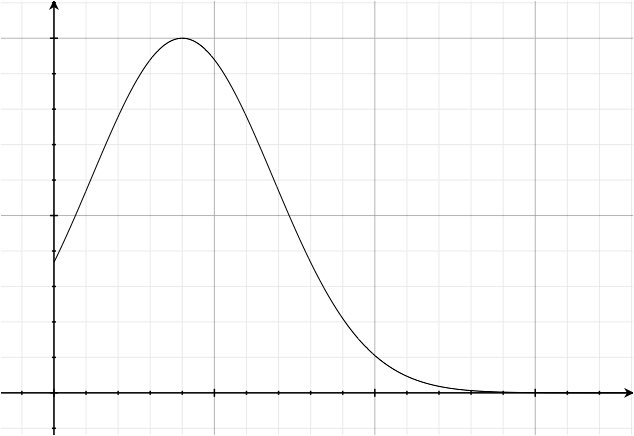}
\caption{} \label{figproc4}
\end{minipage}
\end{figure}

The objective of this paper is to establish a possible geometric approach to the analysis of biological processes and to outline the directions for the future applications of this approach. 
The processes we have in mind are those biochemical processes associated, for example, with food sensitivities 
in humans, cell responses to antigens or to other external stimuli that lead to the release of chemicals such as histamine or cytokines or substance P, and the subsequent signal transmission through the network of nerves in human body. 
There is vast literature that pertains to mathematical modeling of various components of this setup, including, for example, the fundamental work \cite{hodghux}. We mention the paper \cite{tysonnovak} that describes the relationship between the standard cell biology 
perspective and mathematical models based on dynamical systems.   
In the usual setting of control theory (see e.g. the discussion at the beginning of \cite{qiubook}), an example of a system  
that we can have in mind is a car driving on a road that is trying to stay close to the prescribed trajectory, in the presence of disturbances (such as the wind or issues with the road surface). In our setting, where the system 
is the human body, assume that there is a continuous function $A(t)$, $t\ge 0$, that quantifies the amplitude of the cumulative adverse cellular response at time $t$. 
We can state our objective as to change the phase space of our dynamical system, so that with the resulting change of parameters, the response behaviour changes from that on Fig. \ref{figproc1} or Fig. \ref{figproc2} to 
Fig. \ref{figproc4}. 

\begin{figure}[htb]
\centerline{\includegraphics[width=2.0in]{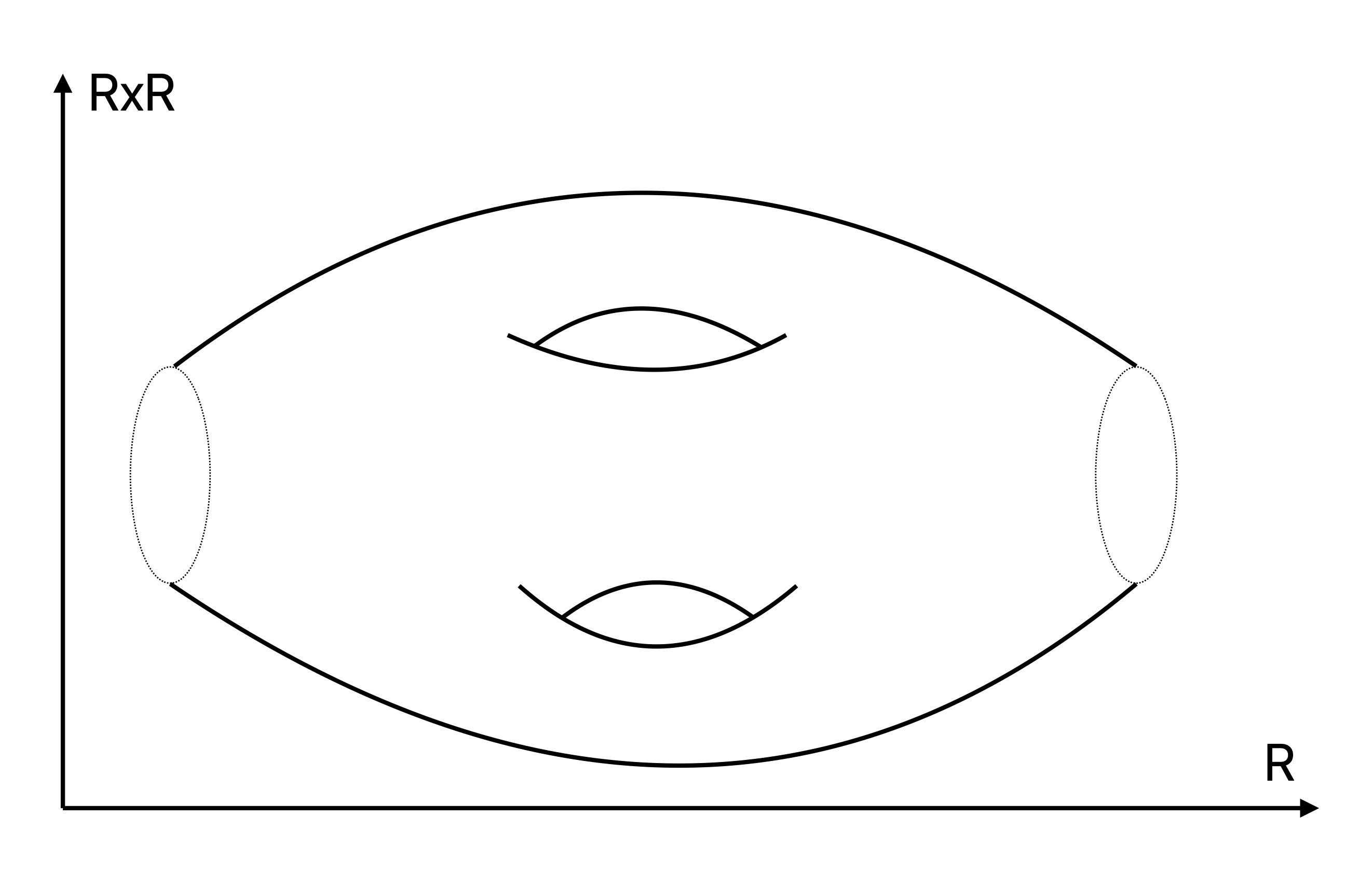}}
\caption{This manifold is not a product of an interval and a $1$-dimensional manifold.} \label{fignotp}
\end{figure}

In all these three Figures, the horizontal axis is time, $t$, the vertical axis is $y$, and the curve is  $y=A(t)$. 
Quantitatively, 
we may say that for a  given $T>0$ and a given $\varepsilon >0$, we would like to have a continuous function $A$ 
with the graph shown in  Fig. \ref{figproc4}, with  
the estimate 
$$
\frac{\int\limits_T^{\infty}|A(t)|dt}{\int\limits_0^T|A(t)|dt}<\varepsilon. 
$$
In this paper, we take the general point of view that a "process" $\Lambda$ is a submanifold of an ambient manifold (which has an extra structure, such as  a Riemannian or pseudoRiemannian metric), as in \cite{barron}, \cite{barronkp}. 
It is not necessarily 
of the form $\Lambda_1\times \Lambda_2$, where $\Lambda_1$ is $1$-dimensional and is supposed to represent "the time". In Fig. \ref{fignotp}, we see an elementary example of a process which is not a direct product. 

In the theory of dynamical systems, the logarithm law,  very informally, states that the time to reach the target 
of size $b>0$ 
in the phase space is asymptotic to $\frac{1}{b}$. See e.g. the discussion and the formal statement in \cite{galatolo}.  
There are many versions of this statement. For example, in \cite{sullivan}, it is proved that in the manifold  $X=H^d/\Gamma $, where $H^d$ is the real hyperbolic space of dimension $d\ge 2$ and $X$ is a finite volume, noncompact quotient by the action of a discrete group of isometries, $\Gamma$, for a point $x_0\in X$ and almost all $v$,  
$$
{\lim\sup}_{t\to\infty}\frac{{\mathrm{dist}} \  v(t)}{t}=\frac{1}{d}
$$
where     ${\mathrm{dist}}  \ v(t)$ denotes the maximum of $1$ and the distance from $x_0$ to the point achieved after traveling time $t$ along the geodesic starting at $x_0$ in direction $v$.  

A revised version of this paper will appear in the Proceedings of the International Conferences on Digital Technology Driven Engineering 2026, Springer. 
  
\section{A process as a manifold}
\label{sec:chemrec}

Approaches that work well in computational biology involve describing complex molecules by positions of its atoms 
in the $3$-dimensional space (as in the protein data bank). This perspective was essential in work of  Baker, Hassabis and Jumper  on the computational protein design and protein structure prediction (Nobel prize in Chemistry 2024 
\cite{nobelp}).    

For our specific purposes, we will associate geometric objects with molecules and processes in a different fashion. 
With a molecule $A$, we associate a compact Riemann surface $\Sigma=\Sigma_A$  of genus $g$, and a holomorphic vector bundle 
\begin{equation}
\label{vectbsigma1}
E=E_{A}= {\mathcal{L}}^{\oplus k}\oplus E_{r_1}\oplus ... \oplus E_{r_m}. 
\end{equation}
Here ${\mathcal{L}}\to\Sigma$ is the canonical bundle, $E_{r_i}$, for $i\in \{ 1,...m\}$, is the trivial vector bundle on $\Sigma$, i.e. $E_{r_i}=\Sigma\times{\mathbb{C}}^{r_i}$. The posiitive integer $k$ is the number of electrons involved in the  covalent or ionic bonds in the molecule. If there are no bonds, i.e. $k=0$, then the summand ${\mathcal{L}}^{\oplus k}$ is not included in (\ref{vectbsigma1}), and instead we have $E=E_{r_1}\oplus ... \oplus E_{r_m}$. 
In the generic case, the positive integer $m$ equals the number of atoms in the molecule. Let's enumerate the atoms: by the index $i$ (an integer), from $i=1$ to $i=m$. The 
positive integer $r_i$, for the $i$-th atom, is defined to be the number of electrons 
in the $i$-th atom that do not participate in any bonds with the other atoms in this molecule.  
 The positive integer $g$ is the sum of atomic numbers $a_i$ of the atoms in the molecule: 
 $$
 g=a_1+...+a_m.
 $$
 We can think that with each individual atom, among these $m$ atoms, there is an associated compact surface 
 $\Sigma^{(i)}$
 of genus $a_i$, and the compact surface $\Sigma$ is obtained by taking connected sums of the surfaces 
 $\Sigma^{(i)}$. 
 The operation of taking a connected sum involves removing a pair of 
 two disks, from two of those component surfaces (each representing one atom), one disk from one surface, and gluing along the boundary circles. This pair of disks represents a bond. The multiplicity of bond, times $2$, will be the contribution of that many copies of 
 the canonical bundle into $E$. 
To rephrase: consider a single covalent bond between two atoms in a molecule. It involves two shared electrons. In this setting it will correspond to removing a small disk from each surface, say removing $D_1$ from $\Sigma^{(i_1)}$, removing $D_2$ from $\Sigma^{(i_2)}$,  
 and gluing   \ $\Sigma^{(i_1)}-D_1$  \ with  \ $\Sigma^{(i_2)}-D_2$, along the  boundary circles. This bond contributes ${\mathcal{L}}^{\oplus 2}$ into $E$. 
 A double covalent bond between two atoms  involves $4$ electrons and corresponds to gluing   \ $\Sigma^{(i_1)}-D_1$  \ with  \ $\Sigma^{(i_2)}-D_2$, along the  boundary circles, with a contribution of ${\mathcal{L}}^{\oplus 4}$ into $E$.
    
 We do not specify the choice of the complex structure on  the Riemann surface $\Sigma$. If it is ever essential that the surfaces $\Sigma^{(i)}$ are Riemann  surfaces, rather than smooth surfaces (i.e. are equipped with complex structure), then   it is reasonable to assume that the choices are made so that the complex structure on $\Sigma$ restricts to the pre-existing complex structure 
 on $\Sigma^{(i)}$ with the disk(s) removed, and the connected sum and gluing operations above extend to gluing 
the canonical bundles on $\Sigma^{(i_1)}-D_1$  and on  \ $\Sigma^{(i_2)}-D_2$ to the canonical bundle 
on the connected sum of $\Sigma^{(i_1)}$  and   \ $\Sigma^{(i_2)}$ (since we have no restriction on the choice of complex structures, and the connected sum is taken finitely many times, it is always possible).  
 
 The surface $\Sigma$ should be a Riemann  surface (a $1$-dimensional complex manifold), rather than a smooth surface (a $2$-dimensional real surface), so we can use the holomorphic line bundles.   

Then, with a chemical reaction of the form
\begin{equation}
\label{chemre1}
\alpha_1A_1+...+\alpha_pA_p \ \to \ \beta_1B_1+...+\beta_qB_q
\end{equation}
where $p,q,\alpha_1,...,\alpha_p,\beta_1,...,\beta_q\in {\mathbb{N}}$, and 
$A_1,...,A_p,B_1,...,B_q$ are molecules, we can associate the vector bundles as follows: 
when $\alpha_1=...=\alpha_p=\beta_1=...=\beta_q=1$, the bundles are  
\begin{equation}
\label{vectbun1}
\pi_1^*E_{A_1}\oplus...\oplus \pi_p^*E_{A_p}\to \Sigma_{A_1}\times ....\times \Sigma _{A_p}
\end{equation}
and 
\begin{equation}
\label{vectbun2}
\pi_1^*E_{B_1}\oplus...\oplus \pi_q^*E_{B_q}\to \Sigma_{B_1}\times ....\times \Sigma _{B_q}
\end{equation}
where for each $j\in \{ 1,..., p\}$, for the molecule $A_j$, 
the vector bundle $E_{A_j}\to \Sigma_{A_j}$ is given by (\ref{vectbsigma1}), \ 
(where $\Sigma$ and $E$ are now  $\Sigma_{A_j}$ and $E_{A_j}$), the map  
$\pi_j$ is the projection onto the $j$-th component of the direct product
$$
\pi_j:\Sigma_{A_1}\times ....\times \Sigma _{A_p}\to \Sigma_{A_j}
$$
$$
(x_1,...,x_p)\mapsto x_j
$$
and $\pi_j^*$ is the pull-back of vector bundles, similarly for $B_1,...,B_q$. Otherwise, for each $\alpha_j>1$, we include 
$\alpha_j$ copies of   $\Sigma_{A_j}$ in the direct product aboive and 
$\alpha _j$ copies of $\pi_j^*E_{A_j}$ in the direct sum above, similarly for $\beta_j$s.  

We can represent the reaction (\ref{chemre1}) as a deformation of (\ref{vectbun1}) into (\ref{vectbun2}) in an appropriate ambient space (e.g. a Riemannian manifold). 

Similar approach can be applied to a biological process that is not, strictly speaking, a chemical reaction: see Section \ref{sec:prote} below.    

\subsection{Methane and chlorine reaction}
\label{sec:metc}
\begin{figure}[htb]
\centering
\begin{minipage}{0.42\textwidth}
\centering
\includegraphics[width=2.2in]{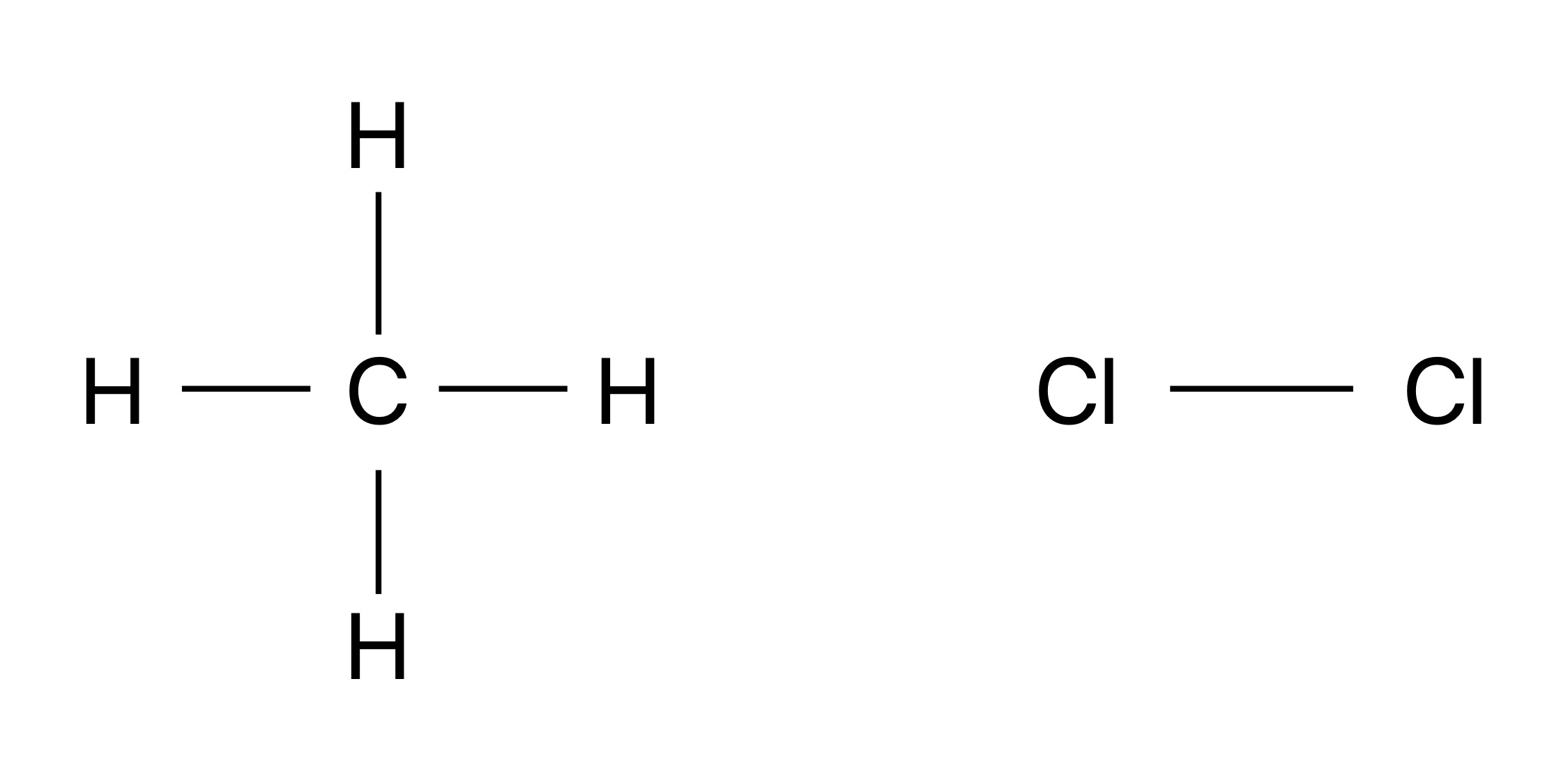}
\caption{$CH_4+Cl_2$} \label{figchem1}
\end{minipage}
\begin{minipage}{0.42\textwidth}
\centering
\includegraphics[width=2.2in]{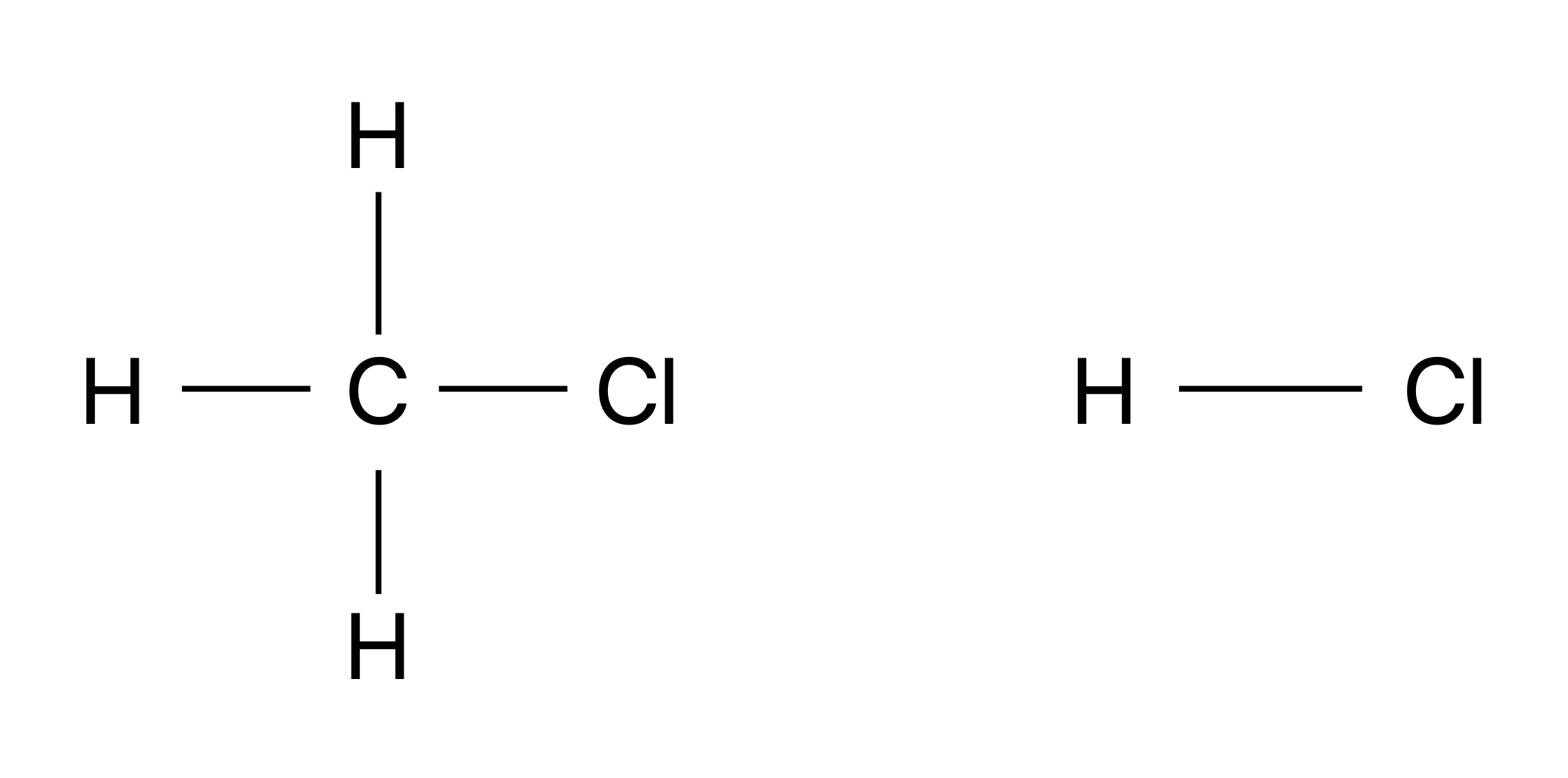}
\caption{$CH_3Cl+HCl$} \label{figchem2}
\end{minipage}
\end{figure}

To illustrate our approach, let's consider the chemical reaction
$$
CH_4+Cl_2 \ \to \ CH_3Cl+HCl.
$$  
We recall the molecular structure: Fig. \ref{figchem1}, \ref{figchem2}. All the bonds are covalent bonds (polar or non-polar). 

Recall the counts for the atoms: 

hydrogen $H$: 1 proton, 1 electron; 

carbon $C$: 6 protons, 6 neutrons, 6 electrons;

chlorine $Cl$: 17 protons, 18 neutrons, 17 electrons. 

The vector bundle that corresponds to  $CH_4+Cl_2$ is 
$$
(\pi_1^* {\mathcal{L}}_{CH_4})^{\oplus 8} \oplus (\pi_2^* {\mathcal{L}}_{Cl_2})^{\oplus 2}
\oplus E_2\oplus E_{32}
\to \Sigma_{CH_4} \times \Sigma_{Cl_2},
$$
where $\Sigma_{CH_4}$ is of genus $10$, $\Sigma_{Cl_2}$ is of genus $34$, 
$E_2$ is the trivial rank $2$ vector bundle on 

$\Sigma_{CH_4} \times \Sigma_{Cl_2}$, and 
$E_{32}\to \Sigma_{CH_4} \times \Sigma_{Cl_2}$ is the trivial vector bundle of rank $32$. 

The vector bundle that corresponds to  $CH_3Cl+HCl$ is 
$$
(\pi_1^* {\mathcal{L}}_{CH_3Cl})^{\oplus 8} \oplus (\pi_2^* {\mathcal{L}}_{HCl})^{\oplus 2}
\oplus E_{18}\oplus E_{16}
\to \Sigma_{CH_3Cl} \times \Sigma_{HCl},
$$
where $\Sigma_{CH_3Cl}$ is of genus $26$, $\Sigma_{HCl}$ is of genus $18$, 
$E_{18}$ is the trivial rank $18$ vector bundle on 
$\Sigma_{CH_3Cl} \times \Sigma_{HCl}$, and 
$E_{16}\to \Sigma_{CH_3Cl} \times \Sigma_{HCl}$ is the trivial vector bundle of rank $16$. 

\subsection{Magnesium oxidation}
\label{sec:magnox}
For a somewhat different example, consider the reaction  
\begin{equation}
\label{reacmo}
2Mg+O_2 \ \to \ 2MgO.
\end{equation}
Recall: 

oxygen $O$: 8 protons, 8 neutrons, 8 electrons; 

magnesium $Mg$: 12 protons, 12 neutrons, 12 electrons. 

The bond in the $O_2$ molecule is a covalent double bond. The bond in the $MgO$ molecule is an ionic bond between 
$Mg^{2+}$ and $O^{2-}$ ions.   

The vector bundle that corresponds to  $Mg+Mg+O_2$ is 
$$
(\pi_3^* {\mathcal{L}}_{O_2})^{\oplus 4}
\oplus E_{12}\oplus E_{12}\oplus E_{12}
\to \Sigma_{Mg} \times \Sigma_{Mg}\times \Sigma_{O_2},
$$
where $\Sigma_{Mg}$ is of genus $12$, $\Sigma_{O_2}$ is of genus $16$, and 
$E_{12}$ is the trivial rank $12$ vector bundle on $\Sigma_{Mg} \times \Sigma_{Mg}\times \Sigma_{O_2}$. 

The vector bundle that corresponds to  $MgO+MgO$ is 
\begin{equation}
\label{vectb2mo}
(\pi_1^* {\mathcal{L}}_{MgO})^{\oplus 2} \oplus (\pi_2^* {\mathcal{L}}_{MgO})^{\oplus 2}
\oplus E_{18}\oplus E_{18}
\to \Sigma_{MgO} \times \Sigma_{MgO},
\end{equation}
where $\Sigma_{MgO}$ is of genus $20$, and 
$E_{18}$ is the trivial rank $18$ vector bundle on 
$\Sigma_{MgO} \times \Sigma_{MgO}$.

\section{Processes involving proteins or enzymes} 
\label{sec:prote}

The approach described in Section \ref{sec:chemrec} can be applied to processes that are more complicated than the reactions in Sections \ref{sec:metc}, \ref{sec:magnox}. 

Various types of adverse  reactions to chemical substances in humans can be modelled as processes that start with 
a finite number of molecules of the types, say $A_1$,..., $A_p$,  and end with a finite number of molecules of  
the types, say $B_1$, ...., $B_q$. The  intermediate steps of the process may involve catalysis, 
binding to receptors, conformational changes in proteins, gene expression, and various subtle factors in cell signalling that may be unknown, or not well understood, or not obvious how to formalize mathematically.  

\subsection{Lactose with or without lactase}
To apply the framework above, for example, to lactose intolerance, we would aim to look for a map   
that takes lactose+lactase (sugar and enzyme), plus water,  to galactose+glucose+lactase, 
and in the absence of lactase,  this process starts with lactose and water, and ends with hydrogen, carbon dioxide, and methane. Recall that the lactose formula is $C_{12}H_{22}O_{11}$, and each of galactose and glucose is 
$C_6H_{12}O_6$. A lactase molecule consists of  21,340 carbon atoms (C),  6,348 oxygen atoms (O),  
5,552 nitrogen atoms (N), and 60 atoms of sulfur (S) \footnote{The protein data bank entries: Lactase  
https://www.rcsb.org/structure/3OB8; Ara h 2 protein:	https://www.rcsb.org/structure/8DB4; 
IgE protein:	https://www.rcsb.org/structure/1O0V}. 

In the spirit of Section \ref{sec:chemrec}, with a molecule of lactose  $C_{12}H_{22}O_{11}$  
we associate a compact Riemann surface of genus 
$$
72+22+88=182,
$$
with a molecule of water $H_2O$ we associate a compact Riemann surface of genus $10$, 
and a lactase molecule  contributes a compact Riemann surface of genus 
$$
6(21340)+8( 6348)+7(5552)+16(60)=218648.
$$
\subsection{Allergic reactions}
\label{sec:allerg}
The analysis of, say, an allergic reaction to peanuts, can be a mathematical model for a process that starts with 
Ara h 2 and IgE proteins and ends with histamine and certain types of cytokines \cite{peanuta}. In an individual without an allergy, the same map 
(with the same domain and codomain),  
starting with an Ara h 2 protein,  should represent a process that does not lead to a release of substances 
responsible for the adverse reaction.  
A molecule of Ara h 2 protein consists of 8878 C atoms,  2792 O atoms,  
2416 N atoms, and 69 S atoms  ${}^\dag$. An IgE molecule consists of 3116 C atoms,  969 O atoms,  
887 N atoms, and 22 S atoms  ${}^\dag$. We can proceed as above. 
\subsection{Other sensitivities to chemicals} 
It is well known that certain types of headaches (in some individuals, but not in other individuals) are triggered by certain foods, or presence of chemicals in the air, 
or by lack of oxygen in the air \cite{oxyg}.   
It is plausible to set up a mathematical model via an approach similar to the above, 
looking for a map $f:X\to Y$ such that for various inputs $X$, yields $f(X)$ that represents an adverse reaction in the individuals sensitive to that set of external factors, and in the other individuals (that are not susceptible) $f(X)$ that represents no response.      

\section{Conclusions and limitations} 
We presented a geometric framework for modeling of chemical and biological processes. In this framework, we associate with a molecule a holomorphic vector bundle on a compact Riemann surface, and suggest to consider a process as a continuous deformation of one vector bundle into another, where time is the deformation parameter. 
We illustrated our approach with two specific examples (methane-chlorine reaction, magnesium oxidation). We indicated how this approach can be used in biomedical applications, to analyize the mechanisms of chemical sensitivities in humans. The well understood issue of lactose intolerance can be used to test our modeling approach. After that, this approach  can be tried with the less understood problems, such as allergies or other adverse reactions (e.g. headaches induced by certain foods or chemicals). 

In our setting, we represent numerical information in geometric form. For example, instead of recording an atomic number $g$, we define a compact surface of genus $g$ (as a remark for general public:  instead of recording a number $1$, we visualize a donut, instead of recording a number $2$, we visualize a pretzel with $2$ holes, and so on).            
Because we also define vector bundles on (products of) these surfaces, we can keep track of other information  
that comes from chemistry or cell biology, by making further choices as needed (e.g. choosing a section or a metric or a connection). If some part of the information is missing or hidden, then we can try to recover this missing part 
from the information that is available to us. 
For example, if, hypothetically, in Section \ref{sec:magnox}
we were given the reaction (\ref{reacmo}), instead,  in the form 
$$
2Mg+A \ \to \ 2MgO
$$
and asked to determine the unknown $A$, then we could have argued that it is a process that starts with the vector bundle 
$$
(\pi_3^* {\mathcal{L}}_{A})^{\oplus k}
\oplus E_{12}\oplus E_{12}\oplus E_{A}
\to \Sigma_{Mg} \times \Sigma_{Mg}\times \Sigma_{A},
$$
where $\Sigma_{Mg}$ is of genus $12$,  
$E_{12}$ is the trivial rank $12$ vector bundle on $\Sigma_{Mg} \times \Sigma_{Mg}\times \Sigma_{A}$, 
and ends with (\ref{vectb2mo}). Then, we conclude that $k=4$, $\Sigma_A$ is of genus $16$, and the rank of $E_A$ is $12$. It takes an extra step to get to the conclusion that $A$ is $O_2$.  In this specific situation it is easier to count atoms, to get to the conclusion $A=O_2$, to answer the question (instead of going through the geometric setup). However, for other types of processes, as in Section \ref{sec:prote}, the setup becomes too complicated to still obey 
the conservation laws. 
It may happen that the answer to questions of the type "find the unknown $A$" would be not unique. One possible way to deal with this would be consider all solutions $A$, together, whether it is a finite set (in the spirit of 
\cite{glavic}) or an appropriate  moduli space of vector bundles on a compact Riemann surface.   

Then, of course, there is the question of control, which, in the biomedical context becomes a complex question of 
finiding a preventative treatment or managing a condition. The diagram Fig. 1.3 \cite{morris2} for glucose control in diabetics is structurally similar to the processes discussed in Section \ref{sec:prote}. 
The problem stated as "hitting the target" in a phase space of a dynamical system in the Introduction can be transformed  into a problem of output estimation with an error uniformly small over all disturbances \cite{morris}.   
This is regardless of the specific mathematical model for the system. The general approaches and techniques such as 
synchronization (as in \cite{gugli}) or some version of a data-driven predictive control method 
(possibly in the spirit of \cite{spsmith}) can be valuable for building intuition and appropriate framework for the challenging problems with a lot of uncertainty in the biomedical applications.

\end{document}